\documentclass[12pt,reqno]{amsart}
\usepackage{amsmath,amsthm,amsbsy,amssymb, longtable,hhline,
hangcaption,caption2}
\usepackage{pstricks,pst-node,pst-tree}
\evensidemargin
\oddsidemargin \marginparwidth 0.5in
\usepackage{amsmath, amssymb, graphicx, array, verbatim, psfig,
epsfig, xspace, lscape}

\newtheorem{theorem}{Theorem}
\newtheorem{lemma}{Lemma}
\newtheorem{corollary}{Corollary}

\usepackage{footnpag} 

\newcommand{\Area}{\mathrm{Area}}

\newcolumntype{C}{>{$}c<{$}}
\newcolumntype{L}{>{$}l<{$}}
\newcolumntype{R}{>{$}r<{$}}

\def\A{\mathcal A}

\newcommand{\T}{\mathcal{T}}
\newcommand{\B}{\mathcal{B}}
\newcommand{\D}{\mathcal{D}}

\newcommand{\LL}{\mathcal{L}}

\newcommand{\V}{\mathcal{V}}
\newcommand{\M}{\mathcal{M}}
\renewcommand{\P}{\mathcal{P}}

\newcommand{\qq}{{\mathbf q}}
\newcommand{\rr}{{\mathbf r}}
\newcommand{\xx}{{\mathbf x}}
\newcommand{\uu}{{\mathbf u}}

\def\RR{\mathbb{R}}
\def\NN{\mathbb{N}}

\newcommand{\F}{{\mathfrak{F}}}
\newcommand{\FQ}{{\mathfrak{F}^{Q}}}

\newcommand{\TT}{{\mathtt{T}}}

\newcommand{\kk}{{\mathbf{k}}}
\def\d{\mathfrak d}
\def\c{\mathfrak c}

\begin{document}

\author[C. Cobeli and A. Zaharescu]
{Cristian Cobeli and Alexandru Zaharescu}

\title[On the Farey fractions with denominators in arithmetic progression]
{On the Farey fractions with denominators in arithmetic progression}\footnotetext{CC is partially supported by the CERES Programme of the Romanian Ministry of Education and Research, contract 4-147/2004.}


\begin{abstract}
Let $\FQ$ be the set of Farey fractions of order $Q$. 
Given the integers $\d\ge 2$ and $0\le \c \le \d-1$, let $\FQ(\c,\d)$ be
the subset of $\FQ$ of those fractions whose denominators are 
$\equiv \c \pmod \d$, arranged in ascending order. The problem we
address here is to show that as $Q\to\infty$, 
there exists a limit probability measuring the distribution of $s$-tuples
of consecutive denominators of fractions in $\FQ(\c,\d)$.
This shows that the clusters of points
$(q_0/Q,q_1/Q,\dots,q_s/Q)\in[0,1]^{s+1}$, where $q_0,q_1,\dots,q_s$
are consecutive denominators of members of $\FQ$ produce a limit set,
denoted by $\D(\c,\d)$.
The shape and the structure of this set are presented in several
particular cases. 
\end{abstract}

\maketitle

\section{Introduction}
This is a continuation of a series of papers dedicated to the study of
the distribution of neighbor denominators of Farey fractions whose
denominators are in arithmetic progression. Previously, we have
treated in \cite{Odds}, \cite{Evens} the cases of pairs of odd and even
denominators, respectively, while here we deal with tuples of consecutive
denominators of fractions in $\FQ(\c,\d)$, the set of Farey
fractions with denominators $\equiv \c \pmod \d$. 
(Here $\c, \d$ are integers, with $\d\ge 2$ and $0\le \c \le \d-1$.)
The motivation for their study comes from their role played in
different problems of various complexities, varying from
applications in the theory of billiards to questions concerned with the
zeros of Dirichlet L-functions.
Although the present work is mostly self-contained, the reader may refer
to the authors~\cite{Evens} and the references within for a wider
introduction of the context and the treatment of some calculations. 

Two generic neighbor fractions from $\FQ$, the set of Farey
fractions of order $Q$, say  $a'/q'$ and $a''/q''$, have two
intrinsic properties. Firstly, the sum $q'+q''$ is always greater than $Q$ and
secondly, $a''q'-a'q''=1$. None of these two properties is generally
true for consecutive members of $\FQ(\c,\d)$, but we shall see that
they may be recovered as initial instances of some more complex connections.

Given a positive integer $s\ge 1$, our main interest lies on the set of
tuples of neighbor denominators of fractions in $\F^{Q}(\c,\d)$:
  \begin{equation*}
       \D_s^{Q}(\c,\d):=\bigg\{(q^0,q^1,\dots,q^s)
	\colon\ 
        \begin{array}{l} 
		q^0,q^1,\dots,q^s\ 	\text{denominators of}
                 \\ \displaystyle
		\text{consecutive fractions in}  \ \F^{Q}(\c,\d)
        \end{array} \bigg\}.
  \end{equation*}  
(Notice that due to some technical constraints, in our notations, the
dimension is $s+1$ and not $s$.)
In fact our aim is to show that there is a limiting set
$\D_s(\c,\d)$ of the scaled set of points 
$\D_s^{Q}(\c,\d)/Q\subset[0,1]^{s+1}$, as $Q\to \infty$.
Strictly speaking, this is the set of limit points of sequences 
$\{\xx_Q\}_{Q\ge 2}$, where each $\xx_Q$ is picked from 
$\D_s^{Q}(\c,\d)/Q$.

More in depth information on $\D_s(\c,\d)$ is revealed if one knows
the concentration of points across its expanse. The answer is given
by Theorem~\ref{Theorem} below, which shows that there exists a 
local density function on $\D_s(\c,\d)$ and gives an explicit
expression for it. Next, let us see the formal definition.
Let $\xx=(x_0,\dots,x_s)$ be a generic point in $[0,1]^{s+1}$ and
denote by $g_s(\xx)=g_s(\xx;\c,\d)$ the function that gives the local
density of  points $(q_0/Q,q_1/Q,\dots,q_s/Q)$ in the
$s+1$-dimensional unit cube, as $Q\to\infty$, 
where $q_0,q_1,\dots,q_s$ are consecutive denominators of fractions in
$\F^{Q}(\c,\d)$.  At any point $\uu=(u_0,\dots,u_s)\in[0,1]^{s+1}$,
we define $g_s(\uu)$ by 
    \begin{equation}\label{eqdefg}
      g_s(\uu):=\lim_{\eta\to 0}
           \frac{\lim\limits_{Q\to \infty}
             \frac{\#\big(\Box\cap\D_s^{Q}(\c,\d)/Q\big)}{\#\D_s^{Q}(\c,\d)}}
                  {4\eta^2}\,,
    \end{equation}
where $\Box\subset\RR^{s+1}$ are cubes  of edge $2\eta$ centered at
$\uu$.

A consequence of the fact that any sequence of consecutive
denominators in $\F^Q$ is uniquely determined by its first two terms
is the framework of $\D_s(\c,\d)$, built as a union of two dimensional
compact surfaces in $\RR^{s+1}$.  This is the reason for which we have
divided in \eqref{eqdefg} by the area of a square of edge $2\eta$
only, and not by $(2\eta)^{s+1}$. 

Thus, in reality $g_s(\uu)=g_s(u,v)$ is a
function of two variables. Here $(u,v)$  will run over a domain that
embodies the Farey series, the Farey triangle with vertices $(0,1)$;
$(1,0)$; $(1,1)$ denoted by $\T$. 

Suppose now that $s=1$, that is, we are in dimension two.
We conclude this introduction with some remarks on the shape of
$\D_1(\c,\d)$, for different $\c$ and $\d$ that we have tested
(see details, tables and pictures in Section~\ref{Annex}).  
It is likely that our observations extend over all $\c$ and
$\d\ge 2$.

The first thing to be remarked is the fact that for any $\d\ge 2$,
$\D_1(\c,\d)$ is a polygon obtained as a union of some sequences of 
polygons with constant local density on each of them, one of them
always including all the others. But the most noteworthy property
is that each of these constant density polygons has a vertex at
$(1,1)$ and looks like a mosaic composed by polygonal pieces, most of
them being quadrangles. The fact that the mosaics exist is not just an
accidental occurrence; on the contrary,  more and more pieces fit
into mosaics with a larger and larger number of components as $\d$
increases.  The mosaics are either symmetric with respect to the first
diagonal or they appear in pairs, whose components are symmetric to
each other with respect to the first diagonal. 

It is not true, as one would guess from tests with many acceptably
small $\d$'s and different $\c$'s, that the larger mosaic (most likely
equal to $\D_1(\c,\d)$) is always a quadrangle. The first
counter-example is  $\D_1(3,12)$, which is a hexagon (see
Fig.~\ref{mosaic_312_11}). 
The shape of the mosaics is more regular when $\d$ has fewer prime
factors. In particular, for each prime modulus $\d$, the exterior
frame of all the mosaics is the same with that from the
case $\d=2$, $\c=0$ (cf.~\cite[Fig. 5, 6]{Evens}), which in turn was the
same in the case $\d=2$, $\c=1$ (cf.~\cite[Fig. 2]{Odds}). 
Mosaics having exactly the same form appear, for example, also when
$\d=4$, $\c=0$, merely than each of them occurs twice.
As opposed to the prime modulus instance, we have included in
Section~\ref{Annex} the pictures that appear in the case $\d=12$,
$\c=3$. One may appreciate these mosaics for their unexpected shape
and beauty.

\section{Notations and Prerequisites}\label{Section1}

Suppose the integer $Q$ is sufficiently large, but fixed. We also fix
$\d\ge 2$, the modulus, $0\le \c \le \d-1$, the residue class, and an
integer $s\ge 1$ ($s+1$ is the dimension).

Then, we define recursively the following objects.
For $0< x,y \le 1$, let 
$\xx^\LL_j(x,y)=x_j^\LL$ be
given by $x_{-1}^\LL=x$, $x_0^\LL=y$ and 
 $x_j^\LL=k_jx_{j-1}^\LL-x_{j-2}^\LL$, for $j\ge 1$,  
where $k_j=\kk_j(x,y):=\Big[\frac{1+x_{j-2}^\LL}{x_{j-1}^\LL}\Big]$. 
We say that $x,y$ are {\em generators} of $\xx^\LL$, and of $\kk$
also, or that $\xx^\LL$ and $\kk$ are {\em generated} by $x,y$.

In order to get a sequence of consecutive denominators of fractions in
$\FQ$ it suffices to know only the first two of them. Moreover, any two
coprime integers $1\le q',q''\le Q$, with $q'+q''> Q$, appear exactly
once in the sequence of consecutive denominators of fractions in
$\FQ$. Then, the subsequent denominators are obtained as follows. Given
two neighbor denominators $1\le q',q''\le Q$, they are succeeded by 
$q_1^\LL,q_2^\LL,\dots$, where
$\qq^\LL_j(q',q'')=q_j^\LL:=k_jq_{j-1}^\LL-q_{j-2}^\LL$, for $j\ge 1$,
and $k_j=\kk_j(q',q'')=\Big[\frac{Q+q_{j-2}^\LL}{q_{j-1}^\LL}\Big]$. We 
put $q_{-1}^\LL=q'$, $q_{0}^\LL=q''$. 
Notice that this values of $k_j$ coincide with those defined above if
$x=q'/Q$ and $y=q''/Q$. Then, in order to simplify the notation, we write
$\kk(q',q'')$ instead of $\kk(q'/Q, q''/Q)$.
As before, we say that $q',q''$ are {\em generators} of $\qq^\LL$ and
of $\kk$ or that $\qq^\LL$ and $\kk$ are {\em generated} by $q',q''$.

A good way to look at a tuple $\kk=(k_1,\dots,k_n)$ is to think that
it is associated to the whole $(n+2)$-tuple $\qq=(q',q'',q_1,\dots,q_n)$
of consecutive denominators in $\FQ$. 
We remark that the link between $\kk$ and $\qq$ is also made by the relations: 
$k_1=(q'+q_1)/q''$, 
$k_2=(q''+q_2)/q_1$, 
$k_3=(q_1+q_3)/q_2$, 
$k_4=(q_2+q_4)/q_3$, etc.

It is plain that $k_j\ge 1$, for $j\ge 1$. 
Additionally, we need to extend the definition of $\kk$ to the case
$n=0$. Then $\kk$ is empty, that is it has no components, 
and we say that its order is zero.

In general, consecutive fractions in $\FQ(\c,\d)$ are not necessarily
consecutive in $\FQ$, but have intercalated in-between several other
fractions from $\FQ$. We remark also that, in general, consecutive
denominators of fractions in $\FQ(\c,\d)$ are not necessarily
coprime. Let $\rr=(r_1,\dots,r_s)$ be an $s$-tuple of positive
integers, and denote $|\rr|=r_1+\cdots+r_s$. 
We say that $\qq=(q_0,\dots,q_s)$, a tuple of
consecutive denominators of fractions in $\FQ(\c,\d)$, is {\em of type}
$\TT(\rr)$ if $\qq\equiv \c\pmod \d$\footnote{We write 
$\qq\equiv \c \pmod \d$ if all the components of $\qq$ are 
$\equiv \c \pmod \d$.} and there exists $(q',q'')$ a pair of consecutive
denominators in $\FQ$ with $\qq^\LL_{-1}(q',q'')=q_0$, 
$\qq^\LL_{-1+r_1}(q',q'')=q_1,\dots, \qq^\LL_{-1+r_1+\cdots+r_s}(q',q'')=q_s$,
and $\qq^\LL_j(q',q'')\not\equiv c \pmod d$, for 
$j\in \big\{0,\dots,|\rr|-1\big\}\setminus 
	\big\{ r_1-1,\dots,r_1+\cdots+r_s-1\big\}$.
In this case we also say that the tuple 
$\kk(q',q'';|\rr|-1)=(k_1,\dots,k_{|\rr|-1})$, with
$k_j=\kk_j(q',q'')$ is {\em of type} $\TT_\rr$. To select the
components that are $\equiv c\pmod d$, we define the 
\emph{choice application} $F^Q_\rr\colon \NN^2 \to\NN^{s+1}$, with
	\begin{equation*}
		F^Q_\rr(q',q''):=\big(q^\LL_{-1},q^\LL_{-1+r_1},\dots,
			q^\LL_{-1+r_1+\cdots+r_s}\big)\,.
	\end{equation*}
Similarly, for $x,y\in(0,1]$, with $x+y>1$, we also put
	\begin{equation*}
		F_\rr(x,y):=\big(x^\LL_{-1},x^\LL_{-1+r_1},\dots,
			x^\LL_{-1+r_1+\cdots+r_s}\big)\,.
	\end{equation*}

Let $\A_\rr(\c,\d)$
be the set of all $\kk(q',q'';|\rr|-1)$ of type $\rr$, for any $q',q''$. 
We remark that the generators of such a $\kk$ are, in general, not
unique. Then, for any $\kk\in\A_\rr(\c,\d)$, we consider the set of
residues relatively prime to $c$, given by 
	\begin{equation*}
	   \M_{\kk,\rr}(\c,\d)=\big\{1\le e \le \d\colon\ 
		\kk(\c,e;|\rr|-1)=\kk\big\}\,.
	\end{equation*}

For example, suppose $Q=25$, $\c=1$, $\d=5$. One can find in $\F^{25}$
the following series of consecutive fractions:
	\begin{equation*}		
	\dots,\frac{7}{16}, \frac{11}{25}, \frac{4}{9}, \frac{9}{20},
	\frac{5}{11}, \frac{11}{24}, \frac{6}{13}, \frac{7}{15}, \frac{8}{17},
	\frac{9}{19}, \frac{10}{21},\dots
	\end{equation*}
From these only $7/16, 5/11, 10/21$ survive, and are consecutive, in
$\F^{25}(1,5)$. Then, in our terminology,  the tuple of denominators
$(16,11,21)$  is of type $\TT(\rr)$, with $\rr=(4,6)$.
In particular, we see that $r_1-1,\dots,r_s-1$ are, respectively, the
number of denominators of consecutive fractions in $\F^Q$ that are
$\not\equiv \c \pmod \d$ intercalated between the fractions with
denominators that are $\equiv \c \pmod \d$. Also,
$\kk(16,25;9)=(1,5,1,4,1,3,2,2,2)\in\A_\rr(1,5)$ has 
$|\rr|-1=4+6-1=9$ components, and
	\begin{equation*}
		F^Q_\rr(16,25)=\big(16,11,21\big)\,.
	\end{equation*}

\section{Lattice Points in Plane Domains}

Given a set $\Omega\subset\RR^2$ and  integers $0\le a,b<\d$, let 
$N'_{a,b;\d}(\Omega)$ be the number of lattice points in $\Omega$ with
relatively prime coordinates congruent modulo $\d$ to $a,b$,
respectively, that is, 
    \begin{equation*}
	N'_{a,b;\d}(\Omega)=\#\big\{
			(m,n)\in\Omega\colon\
				m\equiv a\pmod \d;\ n\equiv b\pmod \d;\ 
					\gcd(m,n)=1
		\big\}\,.
    \end{equation*}

\begin{lemma}\label{LNO}
	Let $R>0$ and let $\Omega\subset\RR^2$ be a convex set of
diameter $\le R$. Let $\d$ be a positive integer and let $0\le a,b<\d$,
with $\gcd(a,b)=1$. Then
    \begin{equation}\label{eqNQ}
	N'_{a,b;\d}(\Omega)=\frac 6{\pi^2 \d^2}\prod_{p\mid \d}\left(1-\frac
			1{p^2}\right)^{-1}
			\Area(\Omega)+O\big(R\log R\big)\,.
    \end{equation}
\end{lemma} 

The proof follows by a standard argument, as in the proof
of~\cite[Lemma 3.1]{BCZ}. 

As a corollary of Lemma~\ref{LNO}, we obtain an asymptotic formula for
the cardinality of $\F^Q(\c,\d)$. For this we use \eqref{eqNQ} and two more
facts. Firstly, the area of the Farey triangle is $\Area(\T^Q)=Q^2/2$
and secondly, the number of residue classes 
$1\le c\le \d$ that are relatively prime to $\c$ is
$\d\cdot\varphi\big(\gcd(\c,\d)\big)/\gcd(\c,\d)$. 
(Here $\varphi(\cdot)$ is the Euler totient function).
Then, we have
	\begin{equation}\label{eqCardFQcd}
		\#\F^Q(\c,\d) = \frac {3Q^2}{\pi^2 }\cdot
		\frac{\varphi\big(\gcd(\c,\d)\big)}{\d\gcd(\c,\d)}
		\prod_{p\mid \d}
			\left(1-\frac 	1{p^2}\right)^{-1}
				+O\big(\d Q\log Q\big)\,.
	\end{equation}

\section{The Density of Points of type $\TT_\rr$}\label{Density}

We count separately the contribution to $g_s(\xx;\c,\d)$ of points of
the same type. Thus, we denote by  $g_\rr(\xx)=g_\rr(\xx;\c,\d)$,
the local density in the unit cube $[0,1]^{s+1}$ of
the points $(q_0/Q,q_1/Q,\dots,q_s/Q)$ of type $\TT_\rr$, as 
$Q\rightarrow\infty$. At any point $\uu=(u_0,\dots,u_s)\in[0,1]^{s+1}$, 
this local density $g_\rr(\uu)$ is defined by 
    \begin{equation}\label{eqdefgr}
      g_\rr(\uu):=\lim_{\eta\to 0}
           \frac{\lim\limits_{Q\to \infty}
             \frac{\#\big(\Box\cap\D_s^Q(\c,\d)/Q\big)}{\#\D_s^Q(\c,\d)}}
                  {4\eta^2}\,,
    \end{equation}
where $\Box\subset\RR^{s+1}$ are cubes  of edge $2\eta$ centered at $\uu$.
Then, we have
  \begin{equation}\label{eqgg}
    g_s(\uu)=\sum_{\rr} g_\rr(\uu)\,,
  \end{equation}
provided we show that each local density $g_\rr(\uu)$ exists, as
$Q\rightarrow\infty$. 
In the following we find each $g_\rr(\uu)$.

\section{The Witness Set}

Let $\eta>0$ be small and let
$\xx^0=(x_0^0,\dots,x_s^0)\in[0,1]^{s+1}$ be the 
point around which we check the density. We consider the
parallelepiped centered at $\xx^0$ and edge $2\eta$ given by
$\square=\square_\eta(\xx^0)
	=(x_0^0-\eta,x_0^0+\eta)\times\dots\times(x_s^0-\eta,x_s^0+\eta)$.
Then, given $\rr=(r_1,\dots,r_s)$, we need to estimate the cardinality
of 
  \begin{equation*}
     \B^Q_\rr(\c,\d)=\left\{(q',q'')\in \NN^2\colon\ 
        \begin{array}{l}  1\le q',q''\le Q,\ \gcd(q',q'') = 1,\ q'+q''>Q,
                 \\ \displaystyle
        \ \kk(q',q'';|\rr|-1)\in\A_\rr(\c,\d),\ 
              F_\rr^Q(q',q''))\in Q\cdot \Box
        \end{array} \right\}.
  \end{equation*}  
This reduces to an area estimate if we put
  \begin{equation*}
     \Omega^Q_\rr(\c,\d)=\left\{(x,y)\in [1,Q]^2\colon\ 
        \begin{array}{l}  
              x+y>Q,\  \\ \displaystyle
         \kk(x,y;|\rr|-1)\in\A_\rr(\c,\d),\  F_\rr(x,y))\in Q\cdot \Box
        \end{array} 
        \right\}.
  \end{equation*}
Then, by Lemma~\ref{LNO},
	\begin{equation}\label{equnu}
	\begin{split}
		\#\B^Q_\rr(\c,\d)
		&=\sum_{e\in\M_{\kk,\rr}(\c,\d)}
			N'_{\c,e;\d}\big(\Omega^Q_\rr(\c,\d)\big)\\
		&=\frac 6{\pi^2\d^2}Q^2 \prod_{p\mid \d}
				\left(1-\frac 1{p^2}\right)^{-1}
			\sum_{e\in\M_{\kk,\rr}(\c,\d)}
			\Area\big(\Omega_\rr(\c,\d)\big)
				+O\big(\d Q\log Q\big)\,,
	\end{split}
	\end{equation}
where
  \begin{equation*}
     \Omega_\rr(\c,\d)=\left\{(x,y)\in (0,1]^2\colon\ 
        \begin{array}{l}  
              x+y>1,\  \\ \displaystyle
         \kk(x,y;|\rr|-1)\in\A_r(\c,\d),\  F_\rr(x,y)\in  \Box
        \end{array} 
        \right\}.
  \end{equation*}

For any $\kk$, we denote 
	\begin{equation*}
		\T_\kk:=\big\{(x,y)\in (0,1]^2\colon\ 
			x+y,\ \kk(x,y;|\rr|-1)=\kk\big\}
	\end{equation*}
and
	\begin{equation}\label{eqbobul}
		\P_\kk(\eta):=\big\{(x,y)\in (0,1]^2\colon\  
		F_\rr(x,y)\in\Box\big\}\,.
	\end{equation}
Also, we put $\T_0:=\T$, the Farey triangle.

Then, we have
	\begin{equation}\label{eqdoi}
		\begin{split}
		\Area\big(\Omega_\rr(\c,\d)\big)=\sum_{\kk\in\A_\rr(\c,\d)}
		\Area\big(\T_\kk\cap\P_\kk(\eta)\big)\,.
		\end{split}
	\end{equation}

By a compactness argument it follows that only finitely many terms of
the series are non-zero, although $\A_\rr(\c,\d)$ may be infinite.
Next we need to see the shape of $\P_\kk(\eta)$, since we are mainly
interested to know $\Area\big(\P_\kk(\eta)\big)$.
This is the object of the next section.

\section{The index $p_r(\kk)$ and the polygon $\P_\kk(\eta)$}
The integer values $k_j$ defined in Section~\ref{Section1} satisfy
the classical {\em mediant} property of the Farey series. For
instance, if $q',q'',q'''$ are consecutive 
denominators of three fractions in $\FQ$, then $k:=(q'+q''')/q''$ is a
positive integer. Hall and Shiu~\cite{HS} called it the {\em index} of
the fractions with denominators $q',q''$, respectively.

More generally, for a series of indices $k_1,k_2,\dots$, 
we consider a sequence of polynomials $p_j(\cdot)$, defined as
follows. Let $p_{-1}(\cdot)=0$, $p_0(\cdot)=1$, and then, for any
$j\ge 1$, 
    \begin{equation}\label{eqR}
	   p_j(k_1,\dots,k_j)
		=k_jp_{j-1}(k_1,\dots,k_{j-1})-p_{j-2}(k_1,\dots,k_{j-2}).
    \end{equation}
The first polynomials with nonempty argument are:
        \begin{align*}
	p_1(\kk)&=k_1;\\
        p_2(\kk)&=k_1k_2-1;\\
        p_3(\kk)&=k_1k_2k_3-k_1-k_3;\\
        p_4(\kk)&=k_1k_2k_3k_4-k_1k_2-k_1k_4-k_3k_4+1;\\
        p_5(\kk)&=k_1k_2k_3k_4k_5-k_1k_2k_3-k_1k_2k_5
			-k_1k_4k_5-k_3k_4k_5+k_1+k_3+k_5.
        \end{align*}  
Often we write $\kk$, meaning the sequence of indices starting with $k_1$,
but notice that the polynomial of rank $j$ depends only on the first
variables $k_1,\dots,k_j$. In particular, one sees that $p_1(\kk)=k_1$
coincides with the index of Hall and Shiu.
Also, we remark the symmetry property: 
        \begin{equation}\label{eqSym}
          p_j(k_j,\dots,k_1)=p_j(k_1,\dots,k_j),\quad\text{for $j\ge 1$}\,.
        \end{equation}
The role played by these polynomials is revealed by the next relation,
which shows that for any $j\ge -1$, $\xx_j^\LL(x,y)$ is a linear
combination of $x$ and $y$:
  \begin{equation}\label{eqLinComb}
    \begin{split}
        \xx_j^\LL(x,y)=p_j(k_1,\dots,k_j)y-p_{j-1}(k_2,\dots,k_j)x\,.
    \end{split}
  \end{equation}

Turning now to the set $\P_\kk$ defined by~\eqref{eqbobul}, where
$\kk=\kk(x,y;|\rr|-1)$, by \eqref{eqLinComb} we see that this is the
set of  points $(x,y)\in\RR^2$ that satisfy simultaneously the
conditions: 
  \begin{equation}\label{eqStrips}
    \begin{split}
        \begin{cases}
                x^0_0-\eta<x<x^0_0+\eta,\\
		        x^0_1-\eta<p_{r_1-1}(k_1,\dots,
		k_{r_1-1})y-p_{r_1-2}(k_2,\dots,k_{r_1-1})x<x^0_1+\eta,\\ 
        x^0_2-\eta<p_{r_1+r_2-1}(k_1,\dots,k_{r_1+r_2-1})y
		-p_{r_1+r_2-2}(k_2,\dots,k_{r_1+r_2-1})x<x^0_2+\eta,\\ 
	\phantom{x^0_2-\eta}\ \ \vdots\\
        	x^0_s-\eta<p_{|\rr|-1}(k_1,\dots,
		k_{|\rr|-1})y-p_{|\rr|-2}(k_2,\dots,k_{|\rr|-1})x<x^0_s+\eta\,.
        \end{cases}
    \end{split}
  \end{equation}
This shows that $\P_\kk(\eta)$ is the intersection of $s+1$ strips and
for $\eta_1,\eta_2>0$ the sets $\P_\kk(\eta_1)$ and $\P_\kk(\eta_2)$ are
similar, the ratio of similarity being equal to $\eta_1/\eta_2$.
Consequently, it follows that 
  \begin{equation}\label{eqAreaP}
    \begin{split}
        \Area(\P_\kk(\eta))=\eta^2\Area(\P_\kk(1))\,,
    \end{split}
  \end{equation}
and $\Area\big(\P_\kk(1)\big)$ is independent of $\eta$.

In particular, in the case $s=1$, for $\kk=(k_1,\dots,k_{r_1-1})$, the
set $\P_\kk(\eta)$ is a parallelogram of center 
  \begin{equation}\label{eqCenter}
    \begin{split}
        C_\kk=\left(x^0_0,\;
	\frac{p_{r_1-2}(k_2,\dots,k_{r_1-1})}
		{p_{r_1-1}(k_1,\dots,k_{r_1-1})}x^0_0
                +\frac{1}{p_{r_1-1}(k_1,\dots,k_{r_1-1})}x^0_1\right),
    \end{split}
  \end{equation}
and area
  \begin{equation}\label{eqAreaParallelogram}
    \begin{split}
        \Area(\P_\kk(\eta))=\frac{4\eta^2}{p_{r_1-1}(\kk)}\,.
    \end{split}
  \end{equation}

\section{The Density of Points of Type $\TT_\rr$}\label{secZeroUnu}

The variable $\eta>0$ is for now fixed, but eventually will tend to
zero. When $\eta\downarrow 0$, for each $\kk$ the polygons
$\P_\kk(\eta)$ are smaller and smaller and converge toward a point
$C_\kk$, which we call the {\em core} of $\P_\kk(\eta)$. If $s=1$,
$\P_\kk(\eta)$ is a parallelogram and the core of $\P_\kk(\eta)$
coincides with its center given by \eqref{eqCenter}. 

Suppose now that $\kk\in\A_\rr(c,d)$ is fixed.
The size of $\Area\big(\T_\kk\cap\P_\kk(\eta)\big)$ depends on the
position of the core with respect to $\T_\kk$. There are three
cases. 

If  $C_\kk\in\overset{\circ}{\T_\kk},$\footnote{For a polygon
$\P\subset\RR^2$, we denote by $\overset{\circ}{\P}$, $\partial\P$ and
$V(\P)$, the topological interior, the boundary, and the set of
vertices of $\P$, respectively.}  then, for $\eta$ small enough,
$\P_\kk(\eta)\subset\T_\kk$, 
so $\Area\big(\T_\kk\cap\P_\kk(\eta)\big)=\Area(\P_\kk(\eta))$. 
Then, by \eqref{eqAreaP}, we get
  \begin{equation}\label{eqInterior}
        \Area(\T_\kk\cap\P_\kk(\eta))=\eta^2\Area(\P_\kk(1))\,,\quad
	\text{if $C_\kk\in\overset{\circ}{\T_\kk}$.}
  \end{equation}

Suppose now that $C_\kk\in\partial\T_\kk\setminus V(\T_\kk)$. Then
there exists a certain bound $\eta_1$ such that if $\eta<\eta_1$, the
intersections $\B_\kk(\eta)=\T_\kk\cap\P_\kk(\eta)$ are polygons similar
to each other. Let $B_\kk$ be the polygon similar to these ones for
which the variable $\eta$ equals $1$ in all the equations of the
boundaries of the strips from \eqref{eqStrips}, whose intersection is
$\P_\kk(\eta)$. So, the size of $B_\kk$ is independent of $\eta$.
Notice that $\B_\kk$ is generally smaller than $\P_\kk(1)$, and even
smaller than $\T_\kk\cap\P_\kk(1)$.
Then, we have
  \begin{equation}\label{eqBoundary}
        \Area(\T_\kk\cap\P_\kk(\eta))=\eta^2\Area(\B_\kk)\,,\quad
	\text{if $C_\kk\in\partial\T_\kk\setminus V(\T_\kk)$ 
		and $\eta<\eta_1$.}
  \end{equation}

If $C_\kk\in V(\T_\kk)$,  the reasoning from the previous case shows
that there exists $\eta_2>0$, with the property that for $\eta_2<\eta$
the polygons $\V_\kk(\eta)=\T_\kk\cap\P_\kk(\eta)$ are similar to each
other. Then, we denote by $\V_\kk$ the polygon similar to these ones for
which $\eta=1$ in all the equations of the
boundaries of the strips from \eqref{eqStrips}. Let us observe that
the size of $\V_\kk$ is independent of $\eta$, and although we use the
same notation, the polygons $\V_\kk$ are distinct for different
vertices of $\T_\kk$. These yield
  \begin{equation}\label{eqVertices}
        \Area(\T_\kk\cap\P_\kk(\eta))=\eta^2\Area(\V_\kk)\,,\quad
	\text{if $C_\kk\in V(\T_\kk)$ and $\eta<\eta_2$.}
  \end{equation}

Inserting the evaluations from \eqref{eqInterior}, \eqref{eqBoundary}
and \eqref{eqVertices} into \eqref{eqdoi}, for 
$0<\eta<\max(\eta_1,\eta_2)$ it yields 
	\begin{equation}\label{eqtrei}
		\begin{split}
		\Area\big(\Omega_\rr(\c,\d)\big) = &
        \eta^2	\sum_{C_\kk\in\overset{\circ}{\T_\kk}} \Area\big(\P_\kk(1)\big)
             +\eta^2\sum_{C_\kk\in\partial{\T_\kk}\setminus 
			V(\T_\kk)} \Area(\B_\kk)\\
        & +\eta^2\sum_{C_\kk\in V(\T_\kk)} \Area(\V_\kk)	\,.
		\end{split}
	\end{equation}

Since the number of tuples $(q_0,\dots,q_s)$ of consecutive
denominators of fractions in $\F^Q(\c,\d)$ is $\#\F^Q(\c,\d)+O(1)$, making
use of \eqref{equnu} and \eqref{eqCardFQcd}, it follows that
   \begin{equation}\label{eqRaportul}
	\begin{split}
	\int\limits_{\Box_\eta(\xx^0)}g_\rr(\xx)\,d\xx
	= & \iint\limits_{\Box_\eta(\xx^0)\cap F_\rr(\T)}
		g_\rr\big(\xx^\LL(x,y)\big)\,dxdy
        = \lim_{Q\rightarrow \infty} \frac{\#\B^Q_\rr(\c,\d)}{\#\F^Q(\c,\d)} \\
	= & \frac{2}{\varphi(\d)}\sum_{e\in\M_{\kk,\rr}(\c,\d)}
			\Area\big(\Omega_\rr(\c,\d)\big)\,.
	\end{split}
   \end{equation}

By Lebesgue differentiation, combining \eqref{eqRaportul} and
\eqref{eqtrei}, we get the following result. 

\begin{theorem}\label{Theorem1}
Let $\d\ge 2$ and $0\le \c \le \d-1$ be integers. Then, 
for any $\xx^0=(x^0_0,x^0_1,\dots,x^0_s)\in [0,1]^{s+1}$, we have:
     \begin{equation}\label{eqgr} 
	\begin{split}
     g_\rr(\xx^0) =& \frac{1}{2\varphi(\d)}\sum_{e\in\M_{\kk,\rr}(\c,\d)}
            	\sum_{C_\kk\in\overset{\circ}{\T_\kk}}
				\Area\big(\P_\kk(1)\big) 
              +\frac{1}{2\varphi(\d)} \sum_{C_\kk\in\partial{\T_\kk}
			\setminus V(\T_\kk)} \Area(\B_\kk) \\              
        & +\frac{1}{2\varphi(\d)}\sum_{C_\kk\in V(\T_\kk)} \Area(\V_\kk) ,
	\end{split}
     \end{equation}
where the sums run over tuples $\kk\in\A_\rr(\c,\d)$.
\end{theorem}

We remark that in \eqref{eqgr},  the first term is essential,  since
it gives the local density on $[0,1]^{s+1}$, except on a set of area
zero.  

\section{The existence of $g_s(\xx)$ and of $\D_s(\c,\d)$}\label{secCompletion}

Putting together the contribution of points of all types, by
\eqref{eqdefgr} and Theorem~\ref{Theorem1}, we get the main result bellow.

\begin{theorem}\label{Theorem}\label{Theorem2}
Let $\d\ge 2$ and $0\le \c \le \d-1$ be integers. Then, for any
$\xx=(x_0,x_1,\dots,x_s)\in [0,1]^{s+1}$, we have 
     \begin{equation*}
	\begin{split}
     g_s(\xx) =& \frac{1}{2\varphi(\d)}\sum_{\rr}\sum_{\kk\in\A_\rr(\c,\d)}
	\sum_{e\in\M_{\kk,\rr}(\c,\d)}
            	\sum_{C_\kk\in\overset{\circ}{\T_\kk}} 
			\Area\big(\P_\kk(1)\big)\\
             &\qquad\qquad+\frac{1}{2\varphi(\d)}
		\sum_{\rr}\sum_{\kk\in\A_\rr(\c,\d)}
	\sum_{C_\kk\in\partial{\T_\kk}\setminus V(\T_\kk)} \Area(\B_\kk)\\
        & +\frac{1}{2\varphi(\d)}\sum_{\rr}\sum_{\kk\in\A_\rr(\c,\d)}
	\sum_{C_\kk\in V(\T_\kk)} \Area(\V_\kk) .
	\end{split}
     \end{equation*}
\end{theorem}
As a consequence, we obtain as a natural object the support set.

\begin{corollary}\label{Corollary}
There exists a limiting set
$\D_s(\c,\d):=\lim_{Q\to\infty}\D_s(\c,\d)/Q$, as $Q\to\infty$.
\end{corollary}

When $s=1$, the theorem can be stated more precisely using
\eqref{eqAreaParallelogram}.  

\begin{corollary}\label{Corollary1}
Let $\d\ge 2$ and $0\le \c \le \d-1$ be integers. Then, 
for any $(x,y)\in [0,1]^{2}$, we have
     \begin{equation}\label{eqgr1} 
	\begin{split}
     g_1(x,y) =& \frac{2}{\varphi(\d)}\sum_{e\in\M_{\kk,r}(\c,\d)}
            	\sum_{C_\kk\in\overset{\circ}{\T_\kk}}\frac{1}{p_{r-1}(\kk)}
             +\frac{1}{\varphi(\d)}\sum_{C_\kk\in\partial{\T_\kk}
				\setminus V(\T_\kk)} 
			\frac{1}{p_{r-1}(\kk)}\\
        & +\frac{1}{2\varphi(\d)}\sum_{C_\kk\in V(\T_\kk)} \Area(\V_\kk) ,
	\end{split}
     \end{equation}
where the sums run over all $r\ge 1$ and $\kk\in\A_r(\c,\d)$.
\end{corollary}
Corollary~\ref{Corollary1} with  $\d=2$  retrieves the
authors~\cite[Theorem 3]{Evens} as a particular instance.

\section{The Mosaics}\label{TabeleEtc}

The noteworthy thing hidden in the background of Theorem~\ref{Theorem2} 
is the geometry of the arrangements of the domains
$F_\rr(\T_\kk)$, which we call \emph{pieces} or \emph{tiles}.
It is easier to see this in the bidimensional case, $s=1$, 
assumed in what follows. Then the tiles are polygons included
in $[0,1]^2$ and the choice application becomes
	\begin{equation*} 
		F_n(x,y)=\big(x, x^{\LL}_{n}(x,y)\big)\,,\quad
		\text{for } n\ge 1\,.
	\end{equation*}

For any $\kk=(k_1,\dots,k_n)$, we shall call \emph{kernel} the integer
$p_n(\kk)$. Moreover, we say that it is the kernel of the tile
$F_n(\T_\kk)$. Notice that the inverse of the kernel is the contribution
of each $\kk$ to $g_1(x,y)$. The tiles of a given kernel fit into a
few larger polygons, which we call \emph{mosaics}. Their common
feature is that always one of their vertices is at $(1,1)$.
They are either symmetric with respect to the first diagonal or they
appear in pairs, symmetric to each other with respect to the first
diagonal.
Most of them are quadrangles, but their shape may vary a lot with
$\d$, $\c$ and the value of the kernel.

These mosaics behave like successive layers of constant density put
over $[0,1]^2$. Then the local density $g_1(x,y)$ at a given point
$(x,y)\in[0,1]^2$ is the sum of the densities on the mosaics (the
inverse of its kernel) stung by $(x,y)$. The contribution to the sum
is halved if $(x,y)$ touches only an edge of a mosaic, and if
the point touches a vertex of a mosaic, it adds to the sum the density
reduced proportionally with the size of the angle of the mosaic at that
vertex. The number of mosaics that lay over $(x,y)\not =(1,1)$ is
finite, and it is endless if $(x,y)=(1,1)$.

For each given $\c, \d$, the number of the mosaics is unbounded, but
their size has a certain rate of decay as their kernel increases.
In the Appendix we have included the larger mosaics in two moduli:
$\d=5$ and $\d=12$. 

It seems that the set $\D_1(\c,\d)$ is always equal to the first mosaic,
which happens to be the largest. This is known to be true when $\d$ is
small and in the cases $\c=0$ and $\d$ prime~\cite{Small}.
Many other intriguing questions are raised by these objects. Here, we
conclude only by pointing out that each of these mosaics has an
associated tree. In the nodes the tree has the tuples $\kk$ that
define the tiles and the arcs link nodes whose corresponding tiles are
adjacent on the mosaic. The root node corresponds to the tile with a
vertex at $(1,1)$. As an example, in Figure~\ref{tree}, it is the tree
associated to the mosaic $SQ_1[9]$ from Figure~\ref{mosaic_25_91}.


\newpage
\section{Appendix--The Plane-mosaics in the cases\\
 $\c=1,2,3,4;\ \d=5$ and $\c=3;\ \d=12$}\label{Annex}

We have assigned names to the mosaics using the following
conventions. The first letter is either $S$ or $N$, according to
whether the mosaic is or not symmetric with respect to the first
diagonal. Since the non-symmetric ones appear in pairs, symmetric to each
other with respect to the first diagonal, we have included the picture
of only one of them. The next letter or group of letters indicates the
shape of the mosaic. The possible configuration are: triangle (T),
quadrangle (Q), pentagon (P), hexagon (H), octagon (O) or concave
hexagon--V-shape (V). The argument is the tuple $\kk$ that gives the
tile from  the North-East corner.  Finally, the subscript represents
the number of components of $\kk$. 

For example,  the mosaic $NP_3[2,2,3]$ (see Figure~\ref{mosaic_25_73}) 
is a non-symmetric pentagon, whose tile from the N-E corner is the
transformation of $\T_{2,2,3}$ through $F_3(x,y)$, and $SQ_1[6]$
(Figure~\ref{mosaic_312_21}) is a symmetric quadrangle, whose piece
from the N-E corner is the image of $\T_{6}$ through $F_1(x,y)$.  

As an exemplification, in Figure~\ref{tree} we have included merely a
tree associated to a mosaic. There, nodes are the tuples $\kk$ defining
the tiles of $NP_3[2,2,3]$ and the arcs connect $\kk$'s that define
adjacent tiles of the mosaic from Figure~\ref{mosaic_25_73}.

More data on the mosaics are entered in Tables~\ref{Table1} and~\ref{Table2}.
On the first column, one can find the kernel, the number whose inverse
gives the local density on the layer given by that mosaic. The entry
on the third column is the number of tiles arranged in the mosaic,
while on the forth are the orders--the number of components--of
$\kk$'s (the smallest and the largest) that produce the tiles. 
In the last column are the coordinates of the vertices of the mosaic.

In the pictures we have used the same color to indicate the chains
of tiles with $\kk$'s of the same orders. There, always neighbor
chains have orders of $\kk$'s that differ by exactly one. 

For the modulus $\d=5$, the mosaics are the same in any of the cases
$\c=1, 2, 3$ or $4$, but they are different when $\c=0$. When $\d=12$,
the situation is more complex, mainly due to the larger number of
factors of $12$.  Due to arithmetical constraints, there are no mosaics
of kernel $2$ when $\d=5$ and $\c=1, 2, 3$ or $4$.

We remark that $\D_1(\c,5)=SQ_0[\cdot]$, for $\c=1,2,3$ or $4$ and
$\D_1(3,12)=SH_1[3]$. In other words, this says that the limiting set
of pairs of consecutive denominators from $\F^{Q}(\c,\d)$ equals, as
$Q\to\infty$, the largest of the mosaics.  

Finally, we mention that $\D_1(3,12)$ is the first case in which
$\D_1(\c,\d)$ is an hexagon, as for any $\d\le 11$ and $0\le \c\le \d$
the set $\D_1(\c,\d)$ has a quadrangular form.

\newpage
\setlength{\doublerulesep}{1pt}
\scriptsize
\smallskip
\setlongtables
\begin{longtable}{|C|C|C|C|L|}
\caption{The mosaics in the cases $\d=5$, $\c=1,2,3$ or $4$.}\label{Table1}\\ \hline
\mathrm{Kernel} & \mathrm{Name} & \mathrm{No.\,  of\, tiles} &
\mathrm{Orders} &
\mathrm{Vertices\, of\, the \, mosaic }    \\ 
\hhline{|=====|}
\endfirsthead
\multicolumn{5}{l}{\small\sl continued from previous page}\\ \hline
\mathrm{Kernel} & \mathrm{Name} & \mathrm{No.\,  of\, tiles} &
\mathrm{Orders} &
\mathrm{Vertices\, of\, the \, mosaic }    \\ 
\hhline{|=====|}
\endhead
\hline
\multicolumn{5}{r}{\small\sl continued on next page} \\ 
\endfoot
\hline
\endlastfoot
   1 & 	SQ_0[\cdot] 	& 21 	& 0-9 	& (1,1); (0,1); (1/6,1/6); (1,0) \\ \hline
   2 & 	-	 	& - 	& - 	& - \\ \hline
   3 & SQ_1[3] 		& 7 	& 1-5 	& (1,1); (2/7,1); (3/8,3/8); (1,2/7) \\ \hline
   4 & SQ_1[4] 		& 27 	& 1-11 	& (1,1); (3/13,1); (2/7,2/7); (1,3/13) \\ \hline
   5 & SHV_4[2,2,2,2]	& 35 	& 4-14 	& (1,1); (1/6,1); (8/43,23/43); (1/2,1/2); (23/43,8/43); (1,1/6) \\ \hline
   6 & SH_1[6]		& 51 	& 1-11 	& (1,1); (1/5,1); (4/19,14/19); (3/8,3/8); (14/19,4/19); (1,1/5) \\ \hline
   7 & NQ_2[2,4]	& 6 	& 2-6 	& (1,1); (6/11,1); (3/5,2/5); (1,3/8) \\ 
   7 & NQ_2[4,2]	& 6 	& 2-6 	& (1,1); (3/8,1); (2/5,3/5); (1,6/11) \\ 
   7 & NP_3[2,2,3]	& 30 	& 3-12 	& (1,1); (3/13,1); (7/17,7/17); (4/5,1/5); (1,6/31) \\ 
   7 & NP_3[3,2,2]	& 30	& 3-12 	& (1,1); (6/31,1); (1/5,4/5); (7/17,7/17); (1,3/13) \\ \hline
   8 & SQ_1[8]		& 21 	& 1-9 	& (1,1); (7/17,1); (4/9,4/9); (1,7/17) \\ 
   8 & SH_3[2,3,2]	& 36	& 3-13 	& (1,1); (7/37,1); (6/31,26/31); (4/9,4/9); (26/31,6/31); (1,7/37) \\ \hline
   9 & SQ_1[9]		& 33 	& 1-9 	& (1,1); (2/7,1); (9/19,9/19); (1,2/7) \\ 
   9 & NP_4[2,2,2,3]	& 7 	& 4-15 	& (1,1); (8/43,1); (7/37,32/37); (1/3,2/3); (1,5/7) \\ 
   9 & NP_4[3,2,2,2]	& 7	& 4-15 	& (1,1); (5/7,1); (2/3,1/3); (32/37,7/37); (1,8/43) \\ \hline
\cdots&	\cdots		&\cdots	&\cdots	&\cdots\\	
\end{longtable}
\normalsize

\captionstyle{flushleft}
\captionwidth=0.49\linewidth
\begin{figure}[b]
    \begin{minipage}[t]{0.49\linewidth}
      \centering
	\includegraphics*[width=\linewidth]{MOZAICURI_2MOD5/mosaic.11}
	\hangcaption{Kernel=1; $\d=5$, $\c=1,2,3$ or $4$.\newline The mosaic $SQ_0[\cdot]$.}\label{mosaic_25_11}
     \end{minipage}
     \hfill
    \begin{minipage}[t]{0.49\linewidth}
      \centering
	\includegraphics*[width=\linewidth]{MOZAICURI_2MOD5/mosaic.31}
	\hangcaption{Kernel=3; $\d=5$, $\c=1,2,3$ or $4$.\newline The mosaic $SQ_1[3]$.}\label{mosaic_25_31}
     \end{minipage}
\end{figure}

\begin{figure}[t]
    \begin{minipage}[t]{0.49\linewidth}
      \centering
	\includegraphics*[width=\linewidth]{MOZAICURI_2MOD5/mosaic.41}
	\hangcaption{Kernel=4; $\d=5$, $\c=1,2,3$ or $4$.\newline The mosaic $SQ_1[4]$.}\label{mosaic_25_41}
     \end{minipage}
     \hfill
    \begin{minipage}[t]{0.49\textwidth}
      \centering
	\includegraphics*[width=\linewidth]{MOZAICURI_2MOD5/mosaic.51}
	\hangcaption{Kernel=5; $\d=5$, $\c=1,2,3$ or $4$.\newline The mosaic $SHV_4[2,2,2,2]$.}\label{mosaic_25_51}
    \end{minipage}
\end{figure}

\begin{figure}[t]
    \begin{minipage}[t]{0.49\textwidth}
      \centering
	\includegraphics*[width=\linewidth]{MOZAICURI_2MOD5/mosaic.61}
	\hangcaption{Kernel=6; $\d=5$, $\c=1,2,3$ or $4$.\newline The mosaic $SH_1[6]$.}\label{mosaic_25_61}
    \end{minipage}
     \hfill
    \begin{minipage}[t]{0.49\textwidth}
      \centering
	\includegraphics*[width=\linewidth]{MOZAICURI_2MOD5/mosaic.71}
	\hangcaption{Kernel=7; $\d=5$, $\c=1,2,3$ or $4$.\newline The mosaic $NQ_2[2,4]$.}\label{mosaic_25_71}
    \end{minipage}
\end{figure}

\begin{figure}[t]
    \begin{minipage}[t]{0.49\linewidth}
      \centering
	\includegraphics*[width=\linewidth]{MOZAICURI_2MOD5/mosaic.73}
	\hangcaption{Kernel=7; $\d=5$, $\c=1,2,3$ or $4$.\newline The mosaic $NP_3[2,2,3]$.}\label{mosaic_25_73}
     \end{minipage}
     \hfill
    \begin{minipage}[t]{0.49\textwidth}
      \centering
	\includegraphics*[width=\linewidth]{MOZAICURI_2MOD5/mosaic.81}
	\hangcaption{Kernel=8; $\d=5$, $\c=1,2,3$ or $4$.\newline The mosaic $SQ_1[8]$.}\label{mosaic_25_81}
    \end{minipage}
\end{figure}

\begin{figure}[t]
    \begin{minipage}[t]{0.49\linewidth}
      \centering
	\includegraphics*[width=\linewidth]{MOZAICURI_2MOD5/mosaic.82}
	\hangcaption{Kernel=8; $\d=5$, $\c=1,2,3$ or $4$.\newline The mosaic $SH_3[2,3,2]$.}\label{mosaic_25_82}
     \end{minipage}
     \hfill
    \begin{minipage}[t]{0.49\textwidth}
      \centering
	\includegraphics*[width=\linewidth]{MOZAICURI_2MOD5/mosaic.91}
	\hangcaption{Kernel=9; $\d=5$, $\c=1,2,3$ or $4$.\newline The mosaic $SQ_1[9]$.}\label{mosaic_25_91}
    \end{minipage}
\end{figure}

\begin{center}
\begin{figure}[t]
    \begin{minipage}[t]{0.49\linewidth}
      \centering
	\includegraphics*[width=\linewidth]{MOZAICURI_2MOD5/mosaic.92}
	\hangcaption{Kernel=9; $\d=5$, $\c=1,2,3$ or $4$.\newline The mosaic $NP_4[2,2,2,3]$.}\label{mosaic_25_92}
     \end{minipage}
\end{figure}
\end{center}

\landscape
	\begin{figure}[H]
	\caption{The tree of $\kk$'s associated to the mosaic
$NP_3[2,2,3]$	($\c=1,2,3$ or $4$ and $\d=5$).} \label{tree} 
	\begin{pspicture}(0,-5)(22,16)
	\psset{framearc=.2}
	\psset{xunit=1.09cm,yunit=1.5cm}
	\rput(9,10){\rnode{3A}{\psframebox{2 2 3}}}
	\rput(9,9){\rnode{4A}{\psframebox{2 3 1 4}}}
	\rput(6,8){\rnode{5A}{\psframebox{3 1 4 1 4}}}
	\rput(12,8){\rnode{5B}{\psframebox{2 3 1 5 1}}}
	\rput(3,7){\rnode{6A}{\psframebox{1 4 1 4 1 4}}}
	\rput(9,7){\rnode{6B}{\psframebox{3 1 4 1 5 1}}}
	\rput(15,7){\rnode{6C}{\psframebox{2 3 1 6 1 2}}}
	\rput(1,6){\rnode{7A}{\psframebox{1 2 4 1 4 1 4}}}
	\rput(5,6){\rnode{7B}{\psframebox{1 4 1 4 1 5 1}}}
	\rput(9,6){\rnode{7C}{\psframebox{3 1 4 2 1 6 1}}}
	\rput(13,6){\rnode{7D}{\psframebox{2 3 2 1 7 1 2}}}
	\rput(17,6){\rnode{7E}{\psframebox{2 3 1 6 1 3 1}}}
	\rput(1,5){\rnode{8A}{\psframebox{1 3 1 5 1 4 1 4}}}
	\rput(7,5){\rnode{8B}{\psframebox{1 4 1 4 2 1 6 1}}}
	\rput(11,5){\rnode{8C}{\psframebox{3 1 4 2 1 7 1 2}}}
	\rput(15,5){\rnode{8D}{\psframebox{2 3 2 1 7 1 3 1}}}
	\rput(5,4){\rnode{9A}{\psframebox{1 4 1 5 1 3 1 6 1}}}
	\rput(9,4){\rnode{9B}{\psframebox{3 1 5 1 3 1 7 1 2}}}
	\rput(13,4){\rnode{9C}{\psframebox{3 1 4 2 1 7 1 3 1}}}
	\rput(17,4){\rnode{9D}{\psframebox{2 3 2 1 8 1 2 3 1}}}
	\rput(3,3){\rnode{10A}{\psframebox{1 4 2 1 6 1 3 1 6 1}}}
	\rput(7,3){\rnode{10B}{\psframebox{1 4 1 5 1 3 1 7 1 2}}}
	\rput(11,3){\rnode{10C}{\psframebox{3 1 5 1 3 1 7 1 3 1}}}
	\rput(15,3){\rnode{10D}{\psframebox{3 1 4 2 1 8 1 2 3 1}}}
	\rput(19,3){\rnode{10E}{\psframebox{2 3 2 1 8 1 2 3 2 1}}}
	\rput(5,2){\rnode{11A}{\psframebox{1 4 2 1 6 1 3 1 7 1 2}}}
	\rput(9,2){\rnode{11B}{\psframebox{1 4 1 5 1 3 1 7 1 3 1}}}
	\rput(13,2){\rnode{11C}{\psframebox{3 1 5 1 3 1 8 1 2 3 1}}}
	\rput(17,2){\rnode{11D}{\psframebox{3 1 4 2 1 8 1 2 3 2 1}}}
	\rput(15,1){\rnode{12A}{\psframebox{3 1 5 1 3 1 8 1 2 3 2 1}}}
	\ncline{3A}{4A}		
	\ncline{4A}{5A}		\ncline{4A}{5B}
	\ncline{5A}{6A}		\ncline{5A}{6B}
	\ncline{5B}{6B}		\ncline{5B}{6C}
	\ncline{6A}{7A}		\ncline{6A}{7B}
	\ncline{6B}{7B}		\ncline{6B}{7C}
	\ncline{6C}{7D} 	\ncline{6C}{7E} 
	\ncline{7A}{8A}
	\ncline{7B}{8B}
	\ncline{7C}{8B} 	\ncline{7C}{8C} 
	\ncline{7D}{8C} 	\ncline{7D}{8D} 
	\ncline{7E}{8D}
	\ncline{8B}{9A}
	\ncline{8C}{9B} 	\ncline{8C}{9C} 
	\ncline{8D}{9C} 	\ncline{8D}{9D} 
	\ncline{9A}{10A} 	\ncline{9A}{10B} 
	\ncline{9B}{10B} 	\ncline{9B}{10C} 
	\ncline{9C}{10C} 	\ncline{9C}{10D} 
	\ncline{9D}{10D} 	\ncline{9D}{10E} 
	\ncline{10A}{11A} 
	\ncline{10B}{11A}   	\ncline{10B}{11B} 
	\ncline{10C}{11B}	\ncline{10C}{11C} 
	\ncline{10D}{11C}	\ncline{10D}{11D}
	\ncline{10E}{11D}
	\ncline{11C}{12A} 
	\ncline{11D}{12A}
	\end{pspicture}
	\end{figure}
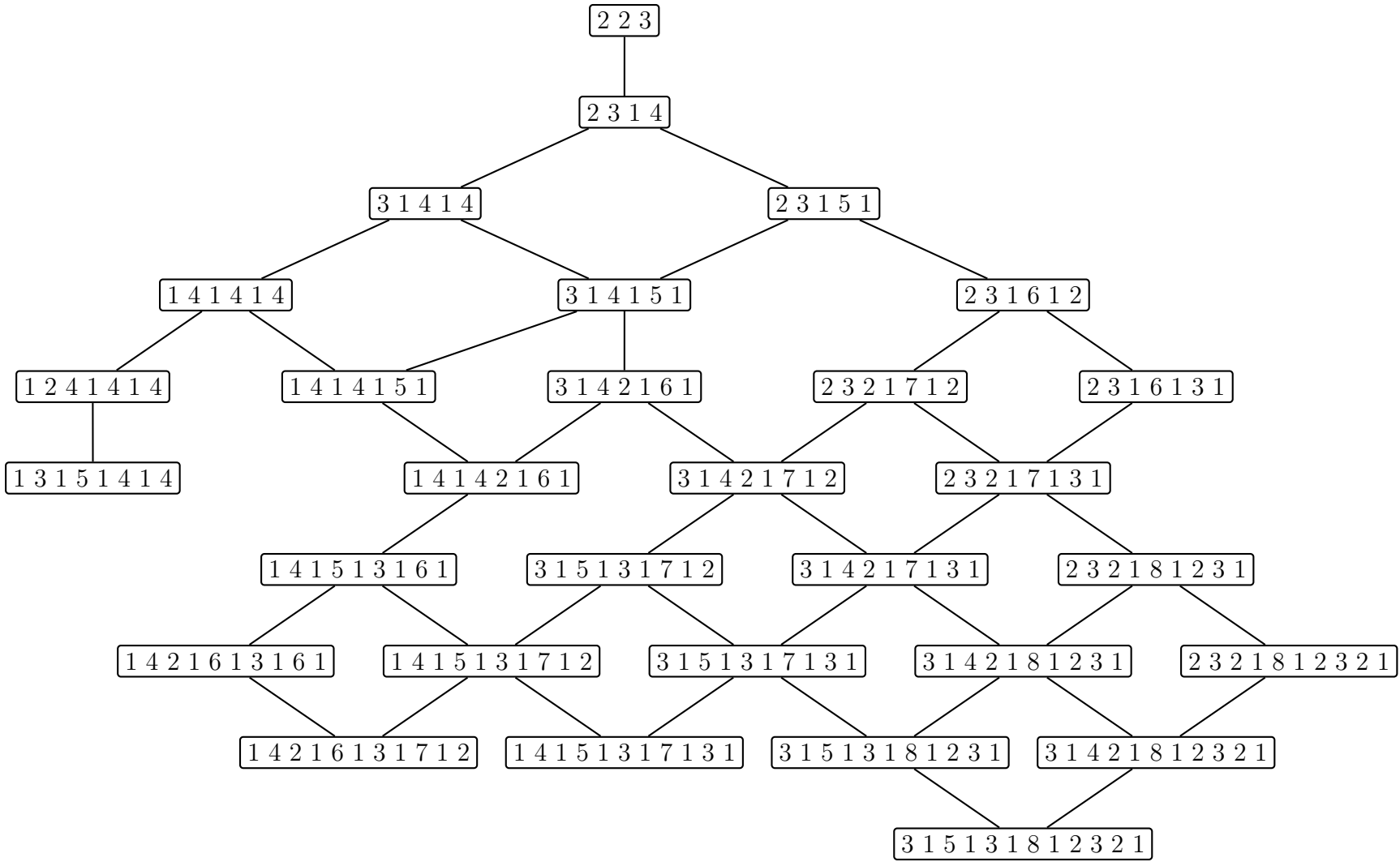
\endlandscape

\landscape
\setlength{\doublerulesep}{1pt}
\scriptsize
\smallskip
\setlongtables
\begin{longtable}{|C|C|C|C|L|}
\caption{The mosaics in the case $\d=12$, $\c=3$.}\label{Table2}\\ \hline
\mathrm{Kernel} & \mathrm{Name} & \mathrm{No.\,  of\, tiles} &
\mathrm{Orders} &
\mathrm{Vertices\, of\, the \, mosaic }    \\ 
\hhline{|=====|}
\endfirsthead
\multicolumn{5}{l}{\small\sl continued from previous page}\\ \hline
\mathrm{Kernel} & \mathrm{Name} & \mathrm{No.\,  of\, tiles} &
\mathrm{Orders} &
\mathrm{Vertices\, of\, the \, mosaic }    \\ 
\hhline{|=====|}
\endhead
\hline
\multicolumn{5}{r}{\small\sl continued on next page} \\ 
\endfoot
\hline
\endlastfoot

   3 & SH_1[3] & \infty & 1-\infty & (1,1); (0,1); (1/13,5/13); (1/5,1/5); (5/13,1/13); (1,0) \\
   3 & SQ_2[2,2] & \infty & 2-\infty & (1,1); (0,1); (1/5,1/5); (1,0) \\ \hline
   6 & SQ_1[6] & 314 & 1-41 & (1,1); (1/9,1); (1/5,1/5); (1,1/9) \\
   6 & SH_5[2,2,2,2,2] & 424 & 5-39 & (1,1); (1/13,1); (1/9,5/9); (1/5,1/5); (5/9,1/9); (1,1/13) \\ \hline
   9 & SQ_1[9] & 63 & 1-13 & (1,1); (1/5,1); (3/7,3/7); (1,1/5) \\
   9 & NQ_4[2,2,2,3] & 56 & 4-18 & (1,1); (1/5,1); (1/4,1/2); (1,1/5) \\
   9 & NQ_4[3,2,2,2] & 56 & 4-18 & (1,1); (1/5,1); (1/2,1/4); (1,1/5) \\ \hline
  12 & SO_5[3,1,6,1,3] & 142 & 5-29 & (1,1); (1/5,1); (3/11,7/11); (4/13, 8/13); (1/3,1/3); (8/13,4/13); (7/11,3/11); (1,1/5) \\
  12 & SQ_{11}[2,\dots,2] & 21 & 11-20 & (1,1); (1/5,1); (1/2,1/2); (1,1/5) \\ \hline
  15 & SQ_1[15] & 38 & 1-11 & (1,1); (3/7,1); (5/9,5/9); (1,3/7) \\
  15 & SH_7[2,3,1,5,1,3,2] & 173 & 7-29 & (1,1); (1/4,1); (3/11,7/11); (5/13,5/13); (7/11,3/11); (1,1/4) \\
  15 & NQ_4[3,2,1,10] & 117 & 4-27 & (1,1); (3/11,1); (5/13,5/13); (1,1/3) \\
  15 & NQ_4[10,1,2,3] & 117 & 4-27 & (1,1); (1/3,1); (5/13,5/13); (1,3/11) \\
  15 & NQ_7[3,2,2,2,2,2,2] & 35 & 7-20 & (1,1); (1/3,1); (7/19,11/19);  (1,1/2) \\
  15 & NQ_7[2,2,2,2,2,2,3] & 35 & 7-20 & (1,1);  (1/2,1); (11/19,7/19); (1,1/3) \\ \hline
  18 & SQ_1[18] & 246 & 1-26 & (1,1); (1/5,1); (3/7,3/7);  (1,1/5) \\
  18 & NQ_6[2,3,1,5,1,4] & 128 & 6-28 & (1,1); (1/5,1); (3/7,3/7);  (1,5/13) \\
  18 & NQ_6[4,1,5,1,3,2] & 128 & 6-28 &(1,1); (5/13,1); (3/7,3/7);  (1,1/5) \\ \hline
  21 & SQ_1[21] & 36 & 1-11 & (1,1); (5/9,1); (7/11,7/11);  (1,5/9) \\
  21 & NT_5[3,2,2,1,14] & 64 & 5-28 & (1,1); (1/3,1); (1,1/5) \\
  21 & NT_5[14,1,2,2,3] & 64 & 5-28 & (1,1);  (1/5,1); (1,1/3) \\
  21 & SQ_7[4,2,1,7,1,2,4] & 16 & 7-13 & (1,1); (5/9,1); (7/11,7/11); (1,5/9) \\
  21 & NQ_8[2,2,3,1,4,2,1,7] & 12 & 8-14 & (1,1); (3/5,1); (7/11,7/11); (1,1/2) \\
  21 & NQ_8[7,1,2,4,1,3,2,2] & 12 & 8-14 & (1,1); (1/2,1); (7/11,7/11); (1,3/5) \\
  21 & NQ_{10}[2,\dots,2,3] & 22 & 10-21 & (1,1); (1/2,1); (5/7,3/7); (1,5/13) \\
  21 & NQ_{10}[3,2,\dots,2] & 22 & 10-21 & (1,1); (5/13,1); (3/7,5/7); (1,1/2) \\ \hline
  24 & SQ_9[6,1,3,1,6,1,3,1,6] & 79 & 9-26 & (1,1); (5/13,1); (1/2,1/2); (1,5/13) \\
  24 & SHV_9[2,2,3,1,5,1,3,2,2] & 161 & 9-33 & (1,1); (1/5,1); (5/13,9/13); (2/3,2/3); (9/13,5/13); (1,1/5)\\
  24 & SH_{11}[4,1,4,1,4,1,4,1,4,1,4] & 94 & 11-27 & (1,1); (7/19,1); (17/43,29/43); (1/2,1/2); (29/43,17/43); (1,7/19) \\ \hline
  27 & SQ_1[27] & 36 & 1-11 & (1,1); (7/11,1); (9/13,9/13); (1,7/11) \\
  27 & NQ_6[10,1,2,3,1,6] & 124 & 6-37 & (1,1); (1/5,1); (11/19,10/19); (1,1/2) \\
  27 & NQ_6[6,1,3,2,1,10] & 124 & 6-37 & (1,1); (1/2,1); (10/19,11/19); (1,1/5) \\
  27 & NQ_8[2,3,2,1,8,1,2,4] & 32 & 8-23 & (1,1); (5/13,1); (7/15,11/15); (1,7/11) \\
  27 & NQ_8[4,2,1,8,1,2,3,2] & 32 & 8-23 & (1,1); (7/11,1); (11/15,7/15); (1,5/13) \\ \hline
\cdots&\cdots&\cdots&\cdots&\cdots\\	
\end{longtable}

\normalsize

\endlandscape
\begin{figure}[t]
    \begin{minipage}[t]{0.49\linewidth}
      \centering
	\includegraphics*[width=\linewidth]{MOZAICURI_3MOD12/mosaic.11}
	\hangcaption{Kernel=3; $\d=12$, $\c=3$.\newline The mosaic $SH_1[3]$.}\label{mosaic_312_11}
     \end{minipage}
     \hfill
    \begin{minipage}[t]{0.49\textwidth}
      \centering
	\includegraphics*[width=\linewidth]{MOZAICURI_3MOD12/mosaic.12}
	\hangcaption{Kernel=3; $\d=12$, $\c=3$.\newline The mosaic $SQ_2[2,2]$.}\label{mosaic_312_12}
    \end{minipage}
\end{figure}

\begin{figure}
    \begin{minipage}[t]{0.49\linewidth}
      \centering
	\includegraphics*[width=\linewidth]{MOZAICURI_3MOD12/mosaic.21}
	\hangcaption{Kernel=6; $\d=12$, $\c=3$.\newline The mosaic $SQ_1[6]$.}\label{mosaic_312_21}
     \end{minipage}
     \hfill
    \begin{minipage}[t]{0.49\textwidth}
      \centering
	\includegraphics*[width=\linewidth]{MOZAICURI_3MOD12/mosaic.22}
	\hangcaption{Kernel=6; $\d=12$, $\c=3$.\newline The mosaic $SQ_5[2,2,2,2,2]$.}\label{mosaic_312_22}
    \end{minipage}
\end{figure}

\begin{figure}
    \begin{minipage}[t]{0.49\linewidth}
      \centering
	\includegraphics*[width=\linewidth]{MOZAICURI_3MOD12/mosaic.31}
	\hangcaption{Kernel=9; $\d=12$, $\c=3$.\newline The mosaic $SQ_1[9]$.}\label{mosaic_312_31}
     \end{minipage}
     \hfill
    \begin{minipage}[t]{0.49\textwidth}
      \centering
	\includegraphics*[width=\linewidth]{MOZAICURI_3MOD12/mosaic.32}
	\hangcaption{Kernel=9; $\d=12$, $\c=3$.\newline The mosaic $NQ_4[2,2,2,3]$.}\label{mosaic_312_32}
    \end{minipage}
\end{figure}

\begin{figure}
    \begin{minipage}[t]{0.49\linewidth}
      \centering
	\includegraphics*[width=\linewidth]{MOZAICURI_3MOD12/mosaic.41}
	\hangcaption{Kernel=12; $\d=12$, $\c=3$.\newline The mosaic $SO_5[3,1,6,1,3]$.}\label{mosaic_312_41}
     \end{minipage}
     \hfill
    \begin{minipage}[t]{0.49\textwidth}
      \centering
	\includegraphics*[width=\linewidth]{MOZAICURI_3MOD12/mosaic.42}
	\hangcaption{Kernel=12; $\d=12$, $\c=3$.\newline The mosaic $SQ_{11}[2,\dots,2]$.}\label{mosaic_312_42}
    \end{minipage}
\end{figure}

\begin{figure}
    \begin{minipage}[t]{0.49\linewidth}
      \centering
	\includegraphics*[width=\linewidth]{MOZAICURI_3MOD12/mosaic.51}
	\hangcaption{Kernel=15; $\d=12$, $\c=3$.\newline The mosaic $SQ_1[15]$.}\label{mosaic_312_51}
     \end{minipage}
     \hfill
    \begin{minipage}[t]{0.49\textwidth}
      \centering
	\includegraphics*[width=\linewidth]{MOZAICURI_3MOD12/mosaic.52}
	\hangcaption{Kernel=15; $\d=12$, $\c=3$.\newline The mosaic $SH_7[2,3,1,5,1,3,2]$.}\label{mosaic_312_52}
    \end{minipage}
\end{figure}

\begin{figure}
    \begin{minipage}[t]{0.49\linewidth}
      \centering
	\includegraphics*[width=\linewidth]{MOZAICURI_3MOD12/mosaic.53}
	\hangcaption{Kernel=15; $\d=12$, $\c=3$.\newline The mosaic $NQ_4[3,2,1,10]$.}\label{mosaic_312_53}
     \end{minipage}
     \hfill
    \begin{minipage}[t]{0.49\textwidth}
      \centering
	\includegraphics*[width=\linewidth]{MOZAICURI_3MOD12/mosaic.55}
	\hangcaption{Kernel=15; $\d=12$, $\c=3$.\newline The mosaic $NQ_7[3,2,2,2,2,2,2]$.}\label{mosaic_312_55}
    \end{minipage}
\end{figure}

\begin{figure}
    \begin{minipage}[t]{0.49\linewidth}
      \centering
	\includegraphics*[width=\linewidth]{MOZAICURI_3MOD12/mosaic.61}
	\hangcaption{Kernel=18; $\d=12$, $\c=3$.\newline The mosaic $SQ_1[18]$.}\label{mosaic_312_61}
     \end{minipage}
     \hfill
    \begin{minipage}[t]{0.49\textwidth}
      \centering
	\includegraphics*[width=\linewidth]{MOZAICURI_3MOD12/mosaic.62}
	\hangcaption{Kernel=18; $\d=12$, $\c=3$.\newline The mosaic $NQ_6[2,3,1,5,1,4]$.}\label{mosaic_312_62}
    \end{minipage}
\end{figure}

\begin{figure}
    \begin{minipage}[t]{0.49\linewidth}
      \centering
	\includegraphics*[width=\linewidth]{MOZAICURI_3MOD12/mosaic.71}
	\hangcaption{Kernel=21; $\d=12$, $\c=3$.\newline The mosaic $SQ_1[21]$.}\label{mosaic_312_71}
     \end{minipage}
     \hfill
    \begin{minipage}[t]{0.49\textwidth}
      \centering
	\includegraphics*[width=\linewidth]{MOZAICURI_3MOD12/mosaic.72}
	\hangcaption{Kernel=21; $\d=12$, $\c=3$.\newline The mosaic $NT_5[3,2,2,1,14]$.}\label{mosaic_312_72}
    \end{minipage}
\end{figure}

\begin{figure}
    \begin{minipage}[t]{0.49\linewidth}
      \centering
	\includegraphics*[width=\linewidth]{MOZAICURI_3MOD12/mosaic.74}
	\hangcaption{Kernel=21; $\d=12$, $\c=3$.\newline The mosaic $SQ_7[4,2,1,7,1,2,4]$.}\label{mosaic_312_74}
     \end{minipage}
     \hfill
    \begin{minipage}[t]{0.49\textwidth}
      \centering
	\includegraphics*[width=\linewidth]{MOZAICURI_3MOD12/mosaic.75}
	\hangcaption{Kernel=21; $\d=12$, $\c=3$.\newline The mosaic $NQ_8[2,2,3,1,4,2,1,7]$.}\label{mosaic_312_75}
    \end{minipage}
\end{figure}

\begin{figure}
    \begin{minipage}[t]{0.49\linewidth}
      \centering
	\includegraphics*[width=\linewidth]{MOZAICURI_3MOD12/mosaic.77}
	\hangcaption{Kernel=21; $\d=12$, $\c=3$.\newline The mosaic $NQ_{10}[2,\dots,2,3]$.}\label{mosaic_312_77}
     \end{minipage}
     \hfill
    \begin{minipage}[t]{0.49\textwidth}
      \centering
	\includegraphics*[width=\linewidth]{MOZAICURI_3MOD12/mosaic.81}
	\hangcaption{Kernel=24; $\d=12$, $\c=3$.\newline $SQ_9[6,1,3,1,6,1,3,1,6]$.}\label{mosaic_312_81}
    \end{minipage}
\end{figure}

\begin{figure}
    \begin{minipage}[t]{0.49\linewidth}
      \centering
	\includegraphics*[width=\linewidth]{MOZAICURI_3MOD12/mosaic.82}
	\hangcaption{Kernel=24; $\d=12$, $\c=3$.\newline  $SHV_9[2,2,3,1,5,1,3,2,2]$.}\label{mosaic_312_82}
     \end{minipage}
     \hfill
    \begin{minipage}[t]{0.49\textwidth}
      \centering
	\includegraphics*[width=\linewidth]{MOZAICURI_3MOD12/mosaic.83}
	\hangcaption{Kernel=24; $\d=12$, $\c=3$.\newline $SH_{11}[4,1,4,1,4,1,4,1,4,1,4]$}\label{mosaic_312_83}
    \end{minipage}
\end{figure}

\begin{figure}
    \begin{minipage}[t]{0.49\linewidth}
      \centering
	\includegraphics*[width=\linewidth]{MOZAICURI_3MOD12/mosaic.91}
	\hangcaption{Kernel=27; $\d=12$, $\c=3$.\newline The mosaic $SQ_1[27]$.}\label{mosaic_312_91}
     \end{minipage}
     \hfill
    \begin{minipage}[t]{0.49\textwidth}
      \centering
	\includegraphics*[width=\linewidth]{MOZAICURI_3MOD12/mosaic.92}
	\hangcaption{Kernel=27; $\d=12$, $\c=3$.\newline The mosaic $NQ_6[10,1,2,3,1,6]$.}\label{mosaic_312_92}
    \end{minipage}
\end{figure}

\begin{figure}
    \begin{minipage}[t]{0.49\linewidth}
      \centering
	\includegraphics*[width=\linewidth]{MOZAICURI_3MOD12/mosaic.94}
	\hangcaption{Kernel=27; $\d=12$, $\c=3$.\newline The mosaic $NQ_8[2,3,2,1,8,1,2,4]$.}\label{mosaic_312_94}
     \end{minipage}
\end{figure}

\end{document}